\documentclass[12pt]{article}
\usepackage[dvips]{graphicx}
\usepackage{amssymb}
\usepackage{amsmath,amsthm}
\usepackage{natbib}

\usepackage{bm}
\bmdefine{\Bzero}{0}
\bmdefine{\Bone}{1}

\def\Bx{{\bf x}}
\def\By{{\bf y}}
\def\Bz{{\bf z}}

\def\Bd{{\bf d}}
\def\Bbeta{{\bf \beta}}

\def\BA{{\rm A}}
\def\BB{{\rm B}}
\def\BC{{\rm C}}
\def\BD{{\rm D}}
\def\BE{{\rm E}}
\def\BF{{\rm F}}
\def\BG{{\rm G}}
\def\BH{{\rm H}}
\def\BI{{\rm I}}
\def\BJ{{\rm J}}
\def\BK{{\rm K}}

\def\Bmu{{\bf \mu}}

\newtheorem{definition}{Definition}[section]

\title{Markov chain Monte Carlo tests for designed experiments}
\author{
Satoshi Aoki\\
Department of Mathematics and Computer Science\\
Kagoshima University\\
and\\
Akimichi Takemura\\
Graduate School of Information Science and Technology\\
University of Tokyo
}
\date{November, 2006}
\begin{document}
\maketitle
\begin{abstract}
We consider conditional exact tests of factor effects in
designed experiments for discrete response variables.  Similarly to
the analysis of contingency tables, a Markov chain Monte Carlo method
can be used for performing exact tests, when large-sample
approximations are poor and the enumeration of the conditional sample
space is
infeasible.  For designed experiments with a single observation for
each run, we formulate log-linear or logistic models and consider
a connected Markov chain over an appropriate sample space.
In particular, we investigate fractional
factorial designs with $2^{p-q}$ runs, noting correspondences to the models
for $2^{p-q}$ contingency tables.
\end{abstract}

\section{Introduction}
Exact calculations of $p$ values for statistical conditional tests
arise mainly in the context of analyzing contingency tables. 
For example, Fisher's exact test is frequently used for evaluating the
hypothesis that the row effect and column effect are independent in the
$2\times 2$ contingency tables. Fisher's exact test is generalized to
$I\times J$ contingency tables in \cite{freeman-halton-1951biometrika}. 
Traditionally, statistical
tests for contingency tables have relied heavily on
large-sample approximations for sampling distribution of the test
statistics. However, many works have shown that large-sample
approximations can be very poor when the contingency table contains
both small and large expected frequencies even when the sample size is
large. See \cite{haberman-1988-jasa}, for example.
Moreover, coupled with rapid development both in computer power and in
techniques of algorithms, exact calculations of $p$ values become
feasible in various settings for
practical use. Consequently, for many types of problems where some ingenious
calculation schemes are invented, 
it is unnecessary to use large-sample approximations for sampling
distributions nowadays when their adequacy is in doubt. 
A typical example is the network algorithm by \cite{mehta-patel-1983jasa} for
calculating exact $p$ values of Freeman-Halton tests in two-way
contingency tables. See the survey paper by \cite{agresti-1992ss}.

At the same time simulation techniques for estimating $p$ values by
Monte Carlo procedures have also developed. In particular, for the
problems where a closed form expression of the sampling distribution can
not be obtained, 
Monte Carlo
method provide powerful tools. Note that, in contrast to the
large-sample approximations, we can estimate $p$ values in arbitrary
accuracy, theoretically, by increasing simulation sizes. 
However, for many models, such as general hierarchical log-linear
models in multi-way contingency tables,  
direct generation of random sample is not straightforward.
In this case, {\it Markov chain
Monte Carlo} techniques can be used.

For performing Markov chain Monte Carlo methods for sampling from
discrete sample space, an important problem is how to construct a
connected Markov chain on the given sample space. Note that, if an
arbitrary connected Markov chain is constructed, the chain can be
modified to give a connected and aperiodic Markov chain with
stationary distribution being the desired null distribution by the
usual Metropolis procedure (\cite{hastings-1970biometrika}, for
example).  As for this
point, the first breakthrough work is given by
\cite{diaconis-sturmfels-1998as}.
The key notion of \cite{diaconis-sturmfels-1998as} is {\it a
  Markov basis}, which enables to construct a connected Markov chain
for arbitrary observed data set.  \cite{diaconis-sturmfels-1998as}
presented a general algorithm for computing a Markov basis in the
settings of a general discrete exponential family of distribution.
Their approach relies on the existence of a Gr\"obner basis of a well
specified polynomial ideal. After \cite{diaconis-sturmfels-1998as}, the
techniques of the Markov chain Monte Carlo method for sampling from
discrete conditional distributions are rapidly developed in the past
decade. See for example \cite{dinwoodie-1998-bernoulli},
\cite{dobra-2003bernoulli}
and the works
by Aoki and Takemura (\cite{aoki-takemura-2003anz},
\cite{aoki-takemura-2003metr},
\cite{aoki-takemura-2005jscs}, \cite{aoki-takemura-2006aism},
\cite{takemura-aoki-2004aism}, \cite{takemura-aoki-2005bernoulli}). 

In this paper, we consider conditional exact tests in the context of
{\it designed experiments} with counts (or ratios of counts)
observations. 
In most of the classical literatures on designed experiments, the
responses are assumed to be normally distributed.  However, in many
practical situations, the experimental data are not normally
distributed. For such non-normal data, the generalized linear models are
frequently used. See \cite{hamada-nelder-1997jqt} or Chapter 13 of
\cite{wu-hamada-2000}, for example. In these literatures, however, exact
testing procedures for non-normal data are not considered. 
Since the experimental design is used when the cost of obtaining the data is
relatively high, 
it is very important to develop techniques of exact 
procedures for the case of non-normal responses. Therefore 
in this manuscript, we consider the exact testing procedures for
non-normal responses, based on the theory of 
the generalized linear models. 
For discrete responses, the above background and strategies also apply,
i.e., to calculate $p$ values for conditional tests,
\begin{enumerate}
\setlength{\itemsep}{0pt}
\item traditionally, large-sample approximations such as the normal
		distribution or
      chi-square distribution are used,
\item if the observed data set contains both small and large expected
      values, the adequacy of the approximation becomes poor,
\item if the sample space is of moderate size, or some ingenious algorithms
      can be used, an exact calculation of $p$ values is possible,
\item if an exact calculation is not feasible, we can rely on Monte
      Carlo procedure,
\item if a closed form expression of the null distribution is not
  given, Markov chain Monte Carlo procedure can be employed, 
  if a Markov basis is available.
\end{enumerate}
The topic we consider in this paper is No.\ 5 of the above
list. In Section 2, we formulate conditional exact tests of factor
effects for fractional factorial designs. 
For designed experiments with a single observation for each
run, we formulate log-linear or logistic models and consider how
to construct a null model to be tested 
using the theory of generalized linear models. 
In Section 3, we consider Markov chain Monte Carlo tests for designed
experiments. First we give a definition of Markov bases and a simple
algorithm for evaluating $p$ values by the Markov chain Monte Carlo
tests in Section 3.1. In Section 3.2, we consider correspondences
between the fractional factorial designs {\it with $2^{p-q}$ runs} and 
models for $2^{p-q}$ contingency tables. We end the paper with some
discussions in Section 4.

\section{Conditional tests for fractional factorial designs}
In this section, we consider the exact conditional tests for the
discrete observations derived from fractional factorial designs. 
We investigate the designs with a single observation for each run, which
is either a count or a ratio of counts. For the former case we consider
the log-linear models, and for the latter case we consider the logistic
models. Since our arguments for the two cases are almost the same, we
first explain in detail the log-linear case in Section 2.1, and then
give only a short
description and a remark for the logistic case in Section 2.2.

\subsection{Exact conditional tests for log-linear models of Poisson observations}
First we investigate the case that the observations are counts of some
events. In this case, it is natural to consider Poisson models. 
To clarify the procedures of exact tests, we take a close look at 
an example of fractional factorial design with counts observations. 
Table 1 is a $1/8$ fraction of a full factorial design (i.e., a
$2^{7-3}$ fractional factorial design) defined from the {\it aliasing relation}
\begin{equation}
\BA\BB\BD\BE = \BA\BC\BD\BF = \BB\BC\BD\BG = \BI,
\label{eqn:aliasing-relation-tab1}
\end{equation}
and response data analyzed in \cite{condra-1993} and
reanalyzed in \cite{hamada-nelder-1997jqt}.
\begin{table*}
\begin{center}
\caption{Design and number of defects $y$ for the wave-solder experiment}
\begin{tabular}{ccccccccrrr}\hline
& \multicolumn{7}{c}{Factor} & \multicolumn{3}{c}{$y$}\\
Run & A & B & C & D & E & F & G & \multicolumn{1}{c}{1} &
 \multicolumn{1}{c}{2} & \multicolumn{1}{c}{3}\\ \hline
 1 & 1 & 1 & 1 & 1 & 1 & 1 & 1 & 13 & 30 &  26\\
 2 & 1 & 1 & 1 & 2 & 2 & 2 & 2 & 4  & 16 &  11\\
 3 & 1 & 1 & 2 & 1 & 1 & 2 & 2 & 20 & 15 &  20\\
 4 & 1 & 1 & 2 & 2 & 2 & 1 & 1 & 42 & 43 &  64\\
 5 & 1 & 2 & 1 & 1 & 2 & 1 & 2 & 14 & 15 &  17\\
 6 & 1 & 2 & 1 & 2 & 1 & 2 & 1 & 10 & 17 &  16\\
 7 & 1 & 2 & 2 & 1 & 2 & 2 & 1 & 36 & 29 &  53\\
 8 & 1 & 2 & 2 & 2 & 1 & 1 & 2 &  5 &  9 &  16\\
 9 & 2 & 1 & 1 & 1 & 2 & 2 & 1 & 29 &  0 &  14\\
10 & 2 & 1 & 1 & 2 & 1 & 1 & 2 & 10 & 26 &   9\\
11 & 2 & 1 & 2 & 1 & 2 & 1 & 2 & 28 &173 &  19\\
12 & 2 & 1 & 2 & 2 & 1 & 2 & 1 &100 &129 & 151\\
13 & 2 & 2 & 1 & 1 & 1 & 2 & 2 & 11 & 15 &  11\\
14 & 2 & 2 & 1 & 2 & 2 & 1 & 1 & 17 &  2 &  17\\
15 & 2 & 2 & 2 & 1 & 1 & 1 & 1 & 53 & 70 &  89\\
16 & 2 & 2 & 2 & 2 & 2 & 2 & 2 & 23 & 22 &   7\\ \hline
\end{tabular}
\end{center}
\end{table*}
In Table 1, the observation $y$ is the number of defects arising in a
wave-soldering process in attaching components to an electronic circuit
card. In Chapter 7 of \cite{condra-1993}, he considered seven factors of a
wave-soldering process: (A) prebake condition, (B) flux density, (C)
conveyer speed, (D) preheat condition, (E) cooling time, (F) ultrasonic
solder agitator and (G) solder temperature, each at two levels with
three boards from each run being assessed for defects.
The aim of this experiment is to decide which levels for each factors
are desirable to reduce solder defects.

In this paper, we only consider designs with a single observation for
each run. Therefore, in this example, we focus on 
the totals for each run in Table 1. This is natural for
the settings of Poisson models, since the set of the totals for each run is the
sufficient statistics for the parameters.
We also ignore the second observation in run 11, which is an obvious outlier as
pointed out in \cite{hamada-nelder-1997jqt}. Therefore the weighted total
of run 11 is $(28 + 19)\times 3/2 = 70.5 \simeq 71$.
By replacing $2$ by $-1$ in Table 1, we
rewrite $k \times p$ design matrix as $D$, where each element is
$+1$ or $-1$. Consequently, we have 
\[
D = \left(\begin{array}{ccccccc}
+1 & +1 & +1 & +1 & +1 & +1 & +1\\
+1 & +1 & +1 & -1 & -1 & -1 & -1\\
+1 & +1 & -1 & +1 & +1 & -1 & -1\\
+1 & +1 & -1 & -1 & -1 & +1 & +1\\
+1 & -1 & +1 & +1 & -1 & +1 & -1\\
+1 & -1 & +1 & -1 & +1 & -1 & +1\\
+1 & -1 & -1 & +1 & -1 & -1 & +1\\
+1 & -1 & -1 & -1 & +1 & +1 & -1\\
-1 & +1 & +1 & +1 & -1 & -1 & +1\\
-1 & +1 & +1 & -1 & +1 & +1 & -1\\
-1 & +1 & -1 & +1 & -1 & +1 & -1\\
-1 & +1 & -1 & -1 & +1 & -1 & +1\\
-1 & -1 & +1 & +1 & +1 & -1 & -1\\
-1 & -1 & +1 & -1 & -1 & +1 & +1\\
-1 & -1 & -1 & +1 & +1 & +1 & +1\\
-1 & -1 & -1 & -1 & -1 & -1 & -1\\
\end{array}
\right),\ \ 
\By = \left(\begin{array}{c}
69\\ 31\\ 55\\ 149\\ 46\\ 43\\ 118\\ 30\\ 
43\\ 45\\ 71\\ 380\\ 37\\ 36\\ 212\\ 52
\end{array}
\right).
\]
Write $\By = (y_1,\ldots,y_k)'$ and $D = (d_{ij}) = (\Bd_1,\ldots,\Bd_p)$
where $\Bd_j = (d_{1j},\ldots,d_{kj})' \in \{-1,+1\}^k$ is the
$j$-th column vector of $D$. $k$ is the number of runs. If there are $q$
aliasing relations defining this design, $k = 2^{p-q}$ holds ($p = 7, q
= 3$ for this example). 
We define $\Bd_{st}$ and $\Bd_{stu}$, $1 \leq s < t < u \leq p$, as
\[
 \Bd_{st} = (d_{1s}d_{1t},\ldots,d_{ks}d_{kt})'
\]
and
\[
 \Bd_{stu} = (d_{1s}d_{1t}d_{1u},\ldots,d_{ks}d_{kt}d_{ku})'
\]
for later use.

The statistical model
for this type of data is constructed from the theory of {\it
generalized linear models} (\cite{mccullagh-nelder-1989}). 
Assume that the observations $y_i$ are mutually independently 
distributed with $\mu_i = E(y_i)$, $i = 1,\ldots,k$. 
The mean parameter  $\mu_i$ is expressed as
\[
g(\mu_i) = \beta_0 + \beta_1x_{i1} + \cdots + \beta_{\nu}x_{i\nu},
\]
where $g(\cdot)$ is the link function and $x_{i1},\ldots,x_{i\nu}$ are
the $\nu$ covariates defined below. 
The sufficient statistic is written as $\sum_{i=1}^{k} x_{ij}y_i$.
The canonical link for
the Poisson distribution is $g(\mu_i) = \log\mu_i$.

Now we define covariates. We write the $\nu$-dimensional
parameter $\beta$ and the covariate matrix $X$ as
\[
 \Bbeta = (\beta_0,\beta_1,\ldots,\beta_{\nu-1})'
\]
and
\[
 X = 
\left(\begin{array}{cccc}
1 & x_{11} & \cdots & x_{1\nu-1}\\
\vdots & \vdots & \cdots & \vdots\\
1 & x_{k1} & \cdots & x_{k\nu-1}
\end{array}
\right) = 
\left(\begin{array}{cccc}
\Bone_k & \Bx_1 & \cdots & \Bx_{\nu-1}
\end{array}
\right) ,
\]
where $\Bone_k = (1,\ldots,1)'$ is the $k$-dimensional column
vector. Since the likelihood function is written as
\[
\begin{array}{ccl}
 \displaystyle\prod_{i=1}^{k}\frac{\mu_i^{y_i}}{y_i!}e^{-\mu_i}
& = & \left(\displaystyle\prod_{i=1}^{k}\frac{e^{-\mu_i}}{y_i!}
\right)\exp\left(\displaystyle\sum_{i=1}^{k}y_i\log\mu_i
\right)\\
& = & \left(\displaystyle\prod_{i=1}^{k}\frac{e^{-\mu_i}}{y_i!}
\right)\exp\left(\beta_0\Bone_k'\By + \displaystyle\sum_{j =
1}^{\nu-1}\beta_{j}\Bx_j'\By\right)\\
& = & \left(\displaystyle\prod_{i=1}^{k}\frac{e^{-\mu_i}}{y_i!}
\right)\exp\left(\Bbeta'X'\By\right),
\end{array}
\]
the sufficient statistic for $\Bbeta$ is $X'\By =
(\Bone_k'\By,\Bx_1'\By,\ldots,\Bx_{\nu-1}'\By)$.

The matrix $X$ is constructed from the design matrix $D$ to reflect
the effects of the factors and their interactions which we intend to
measure.  For example, a simple model which only includes the main
effects for each factor is given as $X = (\Bone_k \ D)$, i.e., $\Bx_j
= \Bd_j$ for $j = 1,\ldots,\nu-1 = p$. On the other hand, we can
consider a more complicated model containing various interaction
effects, under the condition that it is consistent with the 
  aliasing relations.  In this example, the aliasing relation up to
four-factor interactions is derived from
(\ref{eqn:aliasing-relation-tab1}) as follows.
\[
{\small
\begin{array}{l}
\BI = \BA\BB\BD\BE = \BA\BB\BF\BG = \BA\BC\BD\BF = \BA\BC\BE\BG =
 \BB\BC\BD\BG = \BB\BC\BE\BF = \BD\BE\BF\BG\\
\BA = \BB\BD\BE = \BB\BF\BG = \BC\BD\BF = \BC\BE\BG, \qquad 
\BB = \BA\BD\BE = \BA\BF\BG = \BC\BD\BG = \BC\BE\BF\\
\BC = \BA\BD\BF = \BA\BE\BG = \BB\BD\BG = \BB\BE\BF, \qquad 
\BD = \BA\BB\BE = \BA\BC\BF = \BB\BC\BG = \BE\BF\BG\\
\BE = \BA\BB\BD = \BA\BC\BG = \BD\BF\BG, \qquad 
\BF = \BA\BB\BG = \BA\BC\BD = \BB\BC\BF = \BD\BE\BG\\
\BG = \BA\BB\BF = \BA\BC\BE = \BB\BC\BD = \BD\BE\BF\\
\BA\BB = \BD\BE = \BF\BG = \BA\BB\BD\BG = \BA\BC\BD\BG = \BA\BC\BE\BF =
 \BB\BC\BD\BF = \BB\BC\BE\BG\\
\BA\BC = \BD\BF = \BE\BG = \BA\BB\BE\BF = \BB\BC\BD\BE = \BB\BC\BF\BG\\
\BA\BD = \BB\BE = \BC\BF = \BA\BB\BC\BG = \BA\BE\BF\BG = \BB\BD\BF\BG =
 \BC\BD\BE\BG\\
\BB\BC = \BD\BG = \BE\BF = \BA\BB\BD\BF = \BA\BB\BE\BG = \BA\BC\BD\BE =
 \BA\BC\BF\BG\\
\BB\BD = \BA\BE = \BC\BG = \BA\BB\BC\BF = \BA\BD\BF\BG = \BB\BE\BF\BG =
 \BC\BD\BE\BF\\
\BC\BD = \BA\BF = \BB\BG = \BA\BB\BC\BE = \BA\BD\BE\BG = \BB\BD\BE\BF =
 \BC\BE\BF\BG\\
\BA\BG = \BB\BF = \BC\BE = \BA\BB\BC\BD = \BA\BD\BE\BF = \BB\BD\BE\BG =
 \BC\BD\BF\BG\\
\BA\BB\BC = \BA\BD\BG = \BA\BE\BF = \BB\BD\BF = \BB\BE\BG = \BC\BD\BE =
 \BC\BF\BG
\end{array}
}
\]
Subject to the above aliasing relations, 
we may consider appropriate models where all the
parameters are estimable. For example, the saturated model for this
example includes $16 (=k)$ parameters, when $X$ is the Hadamard matrix
of the order $16$. 
One of the interpretations of the saturated
model includes seven main effects, seven
two-factor interaction effects, $\BA\BB, \BA\BC, \BA\BD, \BA\BG, \BB\BC, \BB\BD,
\BC\BD$, and one three-factor interaction effect, $\BA\BB\BC$. We write
this model as
\[\begin{array}{c}
 \BA\BB\BC / \BA\BD / \BB\BD / \BC\BD / \BA\BG / \BE / \BF
\end{array}\]
by the manner of the hierarchical models. 
In this case, as in Table
\ref{tbl:Hadamard-16}, the columns of $X$ can be indexed as
\[
 X = (\Bone_k, \Bd_1,\Bd_2,\Bd_{12},\Bd_3,\Bd_{13},\Bd_{23},\Bd_{123},\Bd_4,\Bd_{14},\Bd_{24},\Bd_5,\Bd_{34},\Bd_6,\Bd_7,\Bd_{17}).
\]
Note that the interpretation of the models 
is not unique. 
An extreme interpretation of the saturated model is
given by ignoring three main effects, $\BE, \BF, \BG$, and considering
all the higher interaction effects among $\BA, \BB, \BC, \BD$, which is
written as
\[
 \BA\BB\BC\BD.
\] 
This interpretation is not realistic in the context of designed
experiments, but it is useful for understanding Markov bases in terms
of contingency tables. 
These two interpretations of the saturated model are shown in 
Table \ref{tbl:Hadamard-16}. 
Note that the concept of identifiable model in the saturated model is 
also given in \cite{pistone-wynn-1996biometrika}. 
Using the theory of Gr\"obner basis, \cite{pistone-wynn-1996biometrika} 
gives a method to obtain one of the maximal identifiable models for a given
design. See also \cite{robbiano-rogantin-1998} and 
\cite{pistone-riccomagno-wynn-2000}.
\begin{table*}
\begin{center}
\caption{Hadamard matrix of the order $16$ and two interpretations.}
\label{tbl:Hadamard-16}
{\scriptsize
\begin{tabular}{|cccccccccccccccc|}\hline
$+1$ & $+1$ & $+1$ & $+1$   & $+1$ & $+1$   & $+1$   & $+1$      & $+1$ & $+1$   & $+1$   & $+1$ & $+1$   & $+1$ & $+1$ & $+1$   \\
$+1$ & $+1$ & $+1$ & $+1$   & $+1$ & $+1$   & $+1$   & $+1$      & $-1$ & $-1$   & $-1$   & $-1$ & $-1$   & $-1$ & $-1$ & $-1$   \\
$+1$ & $+1$ & $+1$ & $+1$   & $-1$ & $-1$   & $-1$   & $-1$      & $+1$ & $+1$   & $+1$   & $+1$ & $-1$   & $-1$ & $-1$ & $-1$   \\
$+1$ & $+1$ & $+1$ & $+1$   & $-1$ & $-1$   & $-1$   & $-1$      & $-1$ & $-1$   & $-1$   & $-1$ & $+1$   & $+1$ & $+1$ & $+1$   \\
$+1$ & $+1$ & $-1$ & $-1$   & $+1$ & $+1$   & $-1$   & $-1$      & $+1$ & $+1$   & $-1$   & $-1$ & $+1$   & $+1$ & $-1$ & $-1$   \\
$+1$ & $+1$ & $-1$ & $-1$   & $+1$ & $+1$   & $-1$   & $-1$      & $-1$ & $-1$   & $+1$   & $+1$ & $-1$   & $-1$ & $+1$ & $+1$   \\
$+1$ & $+1$ & $-1$ & $-1$   & $-1$ & $-1$   & $+1$   & $+1$      & $+1$ & $+1$   & $-1$   & $-1$ & $-1$   & $-1$ & $+1$ & $+1$   \\
$+1$ & $+1$ & $-1$ & $-1$   & $-1$ & $-1$   & $+1$   & $+1$      & $-1$ & $-1$   & $+1$   & $+1$ & $+1$   & $+1$ & $-1$ & $-1$   \\
$+1$ & $-1$ & $+1$ & $-1$   & $+1$ & $-1$   & $+1$   & $-1$      & $+1$ & $-1$   & $+1$   & $-1$ & $+1$   & $-1$ & $+1$ & $-1$   \\
$+1$ & $-1$ & $+1$ & $-1$   & $+1$ & $-1$   & $+1$   & $-1$      & $-1$ & $+1$   & $-1$   & $+1$ & $-1$   & $+1$ & $-1$ & $+1$   \\
$+1$ & $-1$ & $+1$ & $-1$   & $-1$ & $+1$   & $-1$   & $+1$      & $+1$ & $-1$   & $+1$   & $-1$ & $-1$   & $+1$ & $-1$ & $+1$   \\
$+1$ & $-1$ & $+1$ & $-1$   & $-1$ & $+1$   & $-1$   & $+1$      & $-1$ & $+1$   & $-1$   & $+1$ & $+1$   & $-1$ & $+1$ & $-1$   \\
$+1$ & $-1$ & $-1$ & $+1$   & $+1$ & $-1$   & $-1$   & $+1$      & $+1$ & $-1$   & $-1$   & $+1$ & $+1$   & $-1$ & $-1$ & $+1$   \\
$+1$ & $-1$ & $-1$ & $+1$   & $+1$ & $-1$   & $-1$   & $+1$      & $-1$ & $+1$   & $+1$   & $-1$ & $-1$   & $+1$ & $+1$ & $-1$   \\
$+1$ & $-1$ & $-1$ & $+1$   & $-1$ & $+1$   & $+1$   & $-1$      & $+1$ & $-1$   & $-1$   & $+1$ & $-1$   & $+1$ & $+1$ & $-1$   \\
$+1$ & $-1$ & $-1$ & $+1$   & $-1$ & $+1$   & $+1$   & $-1$      & $-1$ & $+1$   & $+1$   & $-1$ & $+1$   & $-1$ & $-1$ & $+1$   \\ \hline
$\BI$& $\BA$& $\BB$&$\BA\BB$& $\BC$&$\BA\BC$&$\BB\BC$&$\BA\BB\BC$& $\BD$&$\BA\BD$&$\BB\BD$& $\BE$&$\BC\BD$& $\BF$& $\BG$&$\BA\BG$\\ 
$\BI$& $\BA$& $\BB$&$\BA\BB$& $\BC$&$\BA\BC$&$\BB\BC$&$\BA\BB\BC$& $\BD$&$\BA\BD$&$\BB\BD$& $\BA\BB\BD$&$\BC\BD$& $\BA\BC\BD$& $\BB\BC\BD$&$\BA\BB\BC\BD$\\ \hline
\end{tabular}
}
\end{center}
\end{table*}

Since the saturated model cannot be tested,  we consider an appropriate
submodel of the saturated model.  For the purpose of illustration we
focus on the model considered in \cite{hamada-nelder-1997jqt}, which is
given as
\begin{equation}
 \BA\BC / \BB\BD / \BE / \BF / \BG, 
\label{eqn:Hamada-Nelder-model}
\end{equation}
i.e., the model of seven main effects and two two-factor interactions. 
We treat this model as the {\it null model} and consider significance tests.
In this case, the null hypothesis can be described as follows.
Permuting the columns of Table \ref{tbl:Hadamard-16}, we 
partition the covariate matrix $X$ of the saturated model as
\[\begin{array}{l}
 X = (X_0\ \ X_1 ),\\
 X_0 = (\Bone_k, \Bd_1,\Bd_2,\Bd_3,\Bd_4,\Bd_5,\Bd_6,\Bd_7,\Bd_{13},\Bd_{24})
    = (\Bone_k, \Bx_1,\ldots,\Bx_{\nu-1}),\\
 X_1 = (\Bd_{12},\Bd_{23},\Bd_{123},\Bd_{14},\Bd_{34},\Bd_{17})
    = (\Bx_{\nu},\ldots,\Bx_{k-1}),\\
\end{array}
\]
and consider the corresponding parameter $\Bbeta =
(\beta_0,\beta_1,\ldots,\beta_{k-1})$.
Then the submodel is specified in the form of a null hypothesis
\[
\mbox{H}_0 :\ \beta_{\nu} = \cdots = \beta_{k-1} = 0.
\]
Under $\mbox{H}_0$, the nuisance parameters are
$\beta_0,\ldots,\beta_{\nu-1}$ and the
sufficient statistic for the nuisance parameters
is $X'_0\By$. 
Then the conditional distribution of $\By$ given the sufficient 
statistics is written as
\begin{equation}
f(\By\ |\ X'_0\By = X'_0\By^o) = C(X'_0\By^o) \displaystyle\prod_{i = 1}^k\frac{1}{y_i!},
\label{eqn:poisson-conditional-distribution}
\end{equation}
where $C(X'_0\By^o)$ is the normalizing constant determined from
$X'_0\By^o$ and written as
\begin{equation}
 C(X'_0\By^o)^{-1} = \displaystyle\sum_{\By \in {\cal F}(X'_0\By^o)}\left(
\displaystyle\prod_{i = 1}^k\frac{1}{y_i!}
\right),
\label{eqn:poisson-constant}
\end{equation}
and
\begin{equation}
 {\cal F}(X'_0\By^o) = \{\By\ |\ X'_0\By = X'_0\By^o,\ y_i \ \mbox{is a
  nonnegative integer for}\ i = 1,\ldots,k\}.
\label{eqn:poisson-fiber}
\end{equation}
Note that by sufficiency the conditional distribution 
does not depend on the values of the nuisance parameters. 

We can now consider significance tests against various alternatives to
$\mbox{H}_0$.  An important alternative is to test the  effect of a 
single additional effect. For example in the
above example we can test presence of 
$\BA\BB$-interaction effect by considering the alternative hypothesis
$\mbox{H}_1: \beta_{\nu}\neq 0$.  Or if we are interested in the
goodness-of-fit test, then the alternative hypothesis is
$\mbox{H}_1: (\beta_{\nu},\dots, \beta_{k-1})\neq (0,\dots,0)$.
Depending on the alternative hypothesis, we can use appropriate test
statistic $T(\By)$, such as the likelihood ratio statistic for testing
$\mbox{H}_0$ against $\mbox{H}_1$.
In Section 3 we give a procedure to
sample from the conditional distribution
(\ref{eqn:poisson-conditional-distribution}).  Therefore we can assess
the conditional distribution of any test statistic $T(\By)$ under 
$\mbox{H}_0$.

For the purpose of illustration we now consider goodness-of-fit
procedures.  Traditional $\chi^2$ tests evaluate the upper probability
for some discrepancy measures such as the deviance, the likelihood
ratio or Pearson goodness-of-fit, based on the asymptotic
distribution, $\chi^2_{k - \nu}$.  For example, the likelihood ratio
statistic
\begin{equation}
\label{eq:g2}
T(\By)= G^2(\By) = 2\sum_{i = 1}^{k}y_i\log\frac{y_i}{\hat{\mu_i}}
\end{equation}
is frequently used, 
where $\hat{\mu_i}$ is the maximum likelihood estimate for $\mu_i$
under the null model (i.e., fitted value), given by
\[
\hat\Bmu= (64.53, 47.25, \dots, 51.42)'
\]
for our example. Then for the observed data $\By^o = (y_1^o,\ldots,y_k^o)'$, 
$T(\By^o)=G^2(\By^o)$ is calculated as $T(\By^o) = 19.096$ and the
corresponding asymptotic $p$ value is $0.0040$ from the asymptotic
distribution $\chi^2_6$. This result tells us that the null hypothesis
is highly significant 
and is rejected. 
Using the conditional distribution
(\ref{eqn:poisson-conditional-distribution}), the exact $p$ value is
written as
\[
 p = \displaystyle\sum_{\By \in {\cal F}(X'_0\By^o)}f(\By\ |\ X'_0\By =
 X'_0\By^o)\Bone(T(\By) \geq T(\By^o)),
\]
where
\[
 \Bone(T(\By) \geq T(\By^o)) = \left\{\begin{array}{ll}
1, & \mbox{if}\ T(\By) \geq T(\By^o),\\
0, & \mbox{otherwise}.
\end{array}
\right.
\]

Unfortunately, however, an enumeration of all the elements in ${\cal
F}(X'_0\By^o)$ and hence the calculation of the normalizing constant 
$C(X'_0\By^o)$ is usually computationally infeasible for large sample
space. Instead, we consider a Markov chain Monte Carlo method described
in Section 3.

\subsection{Exact conditional tests for logistic models of binomial
  observations.}
Next we consider the case that the observation for each run is a ratio of
counts. Table 3 is a $1/2$ fraction of a full factorial design (i.e., a
$2^{4-1}$ fractional factorial design) defined from the relation
\begin{equation}
\BA\BC\BD = \BI
\end{equation}
and response data given by \cite{martin-parker-zenick-1987} and
reanalyzed in \cite{hamada-nelder-1997jqt}. In Table 3, the observation $y$
is the number of good parts out of $1000$ during the stamping process in
manufacturing windshield modeling. The purpose of
\cite{martin-parker-zenick-1987} is to decide the levels for four
factors, (A) poly-film
thickness, (B) oil mixture, (C) gloves and (D) metal blanks, which most
improve the slugging condition. 
\begin{table*}
\begin{center}
\caption{Design and number of good parts $y$ out of $1000$ for the
 windshield molding slugging experiment}
\begin{tabular}{cccccc}\hline
& \multicolumn{4}{c}{Factor} & \\
Run & A & B & C & D & $y$ \\ \hline
1 & 1 & 1 & 1 & 1 & 338\\
2 & 1 & 1 & 2 & 2 & 826\\
3 & 1 & 2 & 1 & 1 & 350\\
4 & 1 & 2 & 2 & 2 & 647\\
5 & 2 & 1 & 1 & 2 & 917\\
6 & 2 & 1 & 2 & 1 & 977\\
7 & 2 & 2 & 1 & 2 & 953\\
8 & 2 & 2 & 2 & 1 & 972\\ \hline
\end{tabular}
\end{center}
\end{table*}
Similarly to Section 2.1, we rewrite this data as
\[
D = \left(\begin{array}{cccc}
+1 & +1 & +1 & +1 \\
+1 & +1 & -1 & -1 \\
+1 & -1 & +1 & +1 \\
+1 & -1 & -1 & -1 \\
-1 & +1 & +1 & -1 \\
-1 & +1 & -1 & +1 \\
-1 & -1 & +1 & -1 \\
-1 & -1 & -1 & +1 \\
\end{array}
\right),\ \ 
\By = \left(\begin{array}{c}
338\\ 826\\ 350\\ 647\\ 917\\ 977\\ 953\\ 972
\end{array}
\right).
\]

As for a statistical model for this type of data, 
it is natural to suppose that the distribution of the observation $y_i$
is the mutually independent binomial distribution ${\rm Bin}(\mu_i, n_i)$,
$i = 1,\ldots,k$, where $n_i = 1000,\ i = 1,\ldots,k = 8$ for this
example. Following the theory of generalized linear models, we also
consider the logit link, which is the canonical link for the binomial
distribution.  It expresses the relation between the mean parameter
$\mu_i$ and the systematic part as
\[
 g(\mu_i) = \log\frac{\mu_i}{1-\mu_i} = \beta_0 + \beta_1x_{1i} + \cdots
 + \beta_{\nu}x_{i\nu}.
\]

The covariate matrix and the corresponding parameters are
defined similarly as Section 2.1, i.e., 
permuting the columns of the Hadamard matrix of the order $8$, 
we write $X = (X_0\ X_1)$ in such a way
that $X_0$ is the covariate matrix for the appropriate null model.
In this example, we consider the
simple main-effects model, $\BA / \BB / \BC / \BD$, which is considered
in \cite{hamada-nelder-1997jqt}. For this model, the covariate matrix is
written as $X_0 = (\Bone_k\ D)$. Similarly to Section 2.1, we consider the
conditional tests for various alternatives to this model. In this
case, since the likelihood function is written as
\[
\begin{array}{ccl}
\displaystyle\prod_{i = 1}^{k}\left(\begin{array}{c}n_i \\
				    y_i\end{array}\right)
 \mu_i^{y_i}(1-\mu_i)^{n_i - y_i} & = & 
\displaystyle\prod_{i = 1}^{k}\left(\begin{array}{c}n_i \\
				    y_i\end{array}\right)
 (1-\mu_i)^{n_i}\left(\frac{\mu_i}{1-\mu_i}\right)^{y_i}\\
& = & \displaystyle\prod_{i = 1}^{k}\left(\begin{array}{c}n_i \\
					  y_i\end{array}\right)
 (1-\mu_i)^{n_i}{\rm exp}(\beta'X_0'\By), 
\end{array}
\]
the nuisance parameters under the null hypothesis are $X_0'\By,
n_1,\ldots,n_k$.
Consequently, the exact conditional tests are based on the conditional
distribution, 
\begin{equation}
f(\By\ |\ X'_0\By = X'_0\By^o, n_1,\ldots,n_k) =
C(X'_0\By^o,n_1,\ldots,n_k) \displaystyle\prod_{i =
1}^k\frac{1}{y_i!(n_i - y_i)!},
\label{eqn:binomial-conditional-distribution}
\end{equation}
where $C(X'_0\By^o, n_1,\ldots,n_k)$ is the normalizing constant
determined from
$X'_0\By^o, n_1,\ldots,n_k$ and written as
\begin{equation}
 C(X'_0\By^o, n_1,\ldots,n_k)^{-1} = \displaystyle\sum_{\By \in {\cal
 F}(X'_0\By^o, n_1,\ldots,n_k)}\left(
\displaystyle\prod_{i = 1}^k\frac{1}{y_i!(n_i - y_i)!}
\right),
\label{eqn:binomial-constant}
\end{equation}
and 
\begin{equation}
 {\cal F}(X'_0\By^o,n_1,\ldots,n_k) = \{\By\ |\ X'_0\By = X'_0\By^o,\
 y_i \in \{0,1,\ldots,n_i\},\ i = 1,\ldots,k\}.
\label{eqn:binomial-fiber}
\end{equation}

For notational convenience, we extend $\By$ to 
\[
 \tilde{\By} = (y_1,\ldots,y_k,n_1-y_1,\ldots,n_k-y_k)'
\]
for the binomial model. Corresponding to this $\tilde{\By}$, we also extend 
$\nu \times k$ matrix $X_0'$ to 
\begin{equation}
\widetilde{X_0}' = \left(\begin{array}{cc}
X_0' & O_{\nu,k}\\
I_k & I_k
\end{array}
\right),
\label{eqn:lawlence-lifting}
\end{equation}
where $O_{\nu,k}$ is the $\nu \times k$ zero matrix and $I_k$ is the
identity matrix of the order $k$. In the theory of the toric ideals, 
$\widetilde{X_0}'$ is called  
the {\it Lawrence lifting} of the
configuration $X_0'$. See \cite{hibi2003} for details. Using $\tilde{\By}$
and $\widetilde{X_0}'$, the condition that $X_0'\By$ and
$n_1,\ldots,n_k$ are fixed is simply written that 
$\widetilde{X_0}'\tilde{\By}$ is fixed. Hereafter
for notational simplicity, we write $\By$ and $X_0'$ instead of
$\tilde{\By}$ and $\widetilde{X_0}'$.
Namely, to express the
conditional distribution and its support for the binomial model, 
we use the expression
(\ref{eqn:poisson-conditional-distribution})(\ref{eqn:poisson-constant})(\ref{eqn:poisson-fiber}),
those for the Poisson model, instead of 
(\ref{eqn:binomial-conditional-distribution})(\ref{eqn:binomial-constant})(\ref{eqn:binomial-fiber}),
respectively,

\section{Markov chain Monte Carlo tests for the designed experiments}
In this section, we consider the Markov chain Monte Carlo methods for
calculating $p$ values defined in Section 2. In Section 3.1, we give an
explanation of Markov chain Monte Carlo methods, along with the
definition of Markov bases. We also describe some algorithms and softwares, 
which are useful in applications. In Section 3.2, we investigate the
relation between the fractional factorial designs and contingency
tables.

\subsection{Markov chain Monte Carlo methods for evaluating $p$ values}
To perform the exact tests defined in Section 2, a useful approach is to
generate samples from the conditional distribution 
$f(\By\ |\ X'_0\By = X'_0\By^o)$ and
calculate $p$ values for any test statistic.
In particular, when the closed form expression of the null distribution
can not be obtained, a Markov chain Monte Carlo approach is a valuable
tool. If a connected Markov chain over 
${\cal F}(X'_0\By^o)$ 
is constructed, the chain can be modified to give a connected and
aperiodic Markov chain with stationary distribution 
$f(\By\ |\ X'_0\By = X'_0\By^o)$ 
by the usual Metropolis procedure. 

To construct a connected chain, a frequently used approach is to
prepare a {\it Markov basis} defined in \cite{diaconis-sturmfels-1998as}. 
Let ${\cal M}(X'_0)$ be the set of integer vectors which are in the
kernel of $X'_0$, i.e.,
\[
{\cal M}(X'_0) = \{\Bz = (z_1,\ldots,z_k)'\ |\ X_0'\Bz = \Bzero, z_i\
\mbox{is an integer for}\ i = 1,\ldots,k\},
\]
where $\Bzero$ is the zero vector. We call the element in ${\cal
M}(X_0')$ 
a {\it move for} $X_0'$, in the sense that adding $\Bz \in
{\cal M}(X'_0)$ to any $\By$ does not change the sufficient statistics,
i.e., 
\[
 X_0'(\By + \Bz) =   X_0'\By.
\]
An important point is that $\By + \Bz$ can include a negative element. 
On the other hand, if $\By + \Bz \in {\cal F}(X'_0\By)$, i.e., $\By + \Bz$ is still a
non-negative vector, we see that $\By$ is {\it moved} to $\By + \Bz \in
{\cal F}(X'_0\By)$ by $\Bz$. Now we give the definition of a Markov
basis.
\begin{definition}
A Markov basis for $X_0'$ is a finite set of moves ${\cal B} =
 \{\Bz_1,\ldots,\Bz_L\}, \Bz_j \in {\cal M}(X_0'), j = 1,\ldots,L$,
such that, for any $\By, \By^* \in {\cal F}(X_0'\By^o)$, 
there exists $A > 0$,
 $(\varepsilon_1,\Bz_{j_1}),\ldots,\allowbreak(\varepsilon_A,\Bz_{j_A})$ with
 $\varepsilon_s \in \{-1,+1\}$, $\Bz_{j_s} \in {\cal B}$, $s =
 1,\ldots,A$, satisfying
\[
\By = \By^* + \sum_{s = 1}^A\varepsilon_s\Bz_{j_s} \ \mbox{and}\ \By^* + 
\sum_{s = 1}^a\varepsilon_s\Bz_{j_s} \in {\cal F}(X_0'\By^o)\
 \mbox{for}\ a = 1,\ldots,A.
\]
\end{definition}

By definition, a Markov basis enables to construct a connected chain
over ${\cal F}(X_0'\By^o)$, which is modified so as to have the null
distribution $f(\By\ |\ X_0'\By = X_0'\By^o)$ 
as the stationary distribution by the Metropolis-Hastings procedure.
Therefore we can perform various conditional tests by the
Monte Carlo sampling. 
We give a simple algorithm to calculate $p$ values for some
test statistic $T(\cdot)$ based on the
Markov chain Monte Carlo sampling.

\bigskip

{\tt
Input:\ Markov basis ${\cal B}$, observed data $\By^o$, covariate matrix
$X_0'$, \\ \hspace*{1.5cm} size of sample $N$, null distribution $f(\cdot)$, 
test statistic $T(\cdot)$

Output:\ $p$ value

Variables:\ obs, count, sig, $\By, \By_{next}$\\

Step 1:\ $\mbox{obs} = T(\By^o)$.\ $\By = \By^o$.\ $\mbox{count} = 0$.\  $\mbox{sig} = 0$.

Step 2:\ Choose $\Bz$ from ${\cal B}$ randomly. $I = \{n\ |\ \By + n\Bz
\in {\cal F}(X_0'\By^o)\}$.


Step 3:\ Select $\By_{next}$ from $\{\By + n\Bz\ |\ n \in I\}$ with
probability
\[
 p_n = \frac{f(\By + n\Bz)}{\displaystyle\sum_{n \in I}f(\By + n\Bz)}.
\]

Step 4:\ If $T(\By_{next}) \geq  \mbox{obs}$ then $\mbox{sig} = \mbox{sig} + 1$.

Step 5:\ $\By = \By_{next}$.\ $\mbox{count} = \mbox{count} + 1$.

Step 6:\ If $\mbox{count} < N$ then Go to Step 2.

Step 7:\ Estimated $p$ value is $\displaystyle\frac{\mbox{sig}}{N}$

}

\bigskip
Note that we need not calculate the normalizing constant, $C(X_0'\By^o)$
in (\ref{eqn:poisson-constant}) or $C(X_0'\By^o,\allowbreak n_1,\ldots,n_k)$ in
(\ref{eqn:binomial-constant}), of the null distribution $f(\cdot)$,
since it is canceled in the numerator and denominator in Step 3. 

Derivation of Markov bases is itself a very interesting problem.
Markov bases can be very complicated and hard to compute for large
models.
Many works, including the original work
by \cite{diaconis-sturmfels-1998as}, have relied on the theory of
computational algebra and Gr\"obner bases. See
\cite{diaconis-sturmfels-1998as}, \cite{dinwoodie-1998-bernoulli},
\cite{dobra-2003bernoulli}. On the
other hand, a series of
works by Aoki and Takemura investigates the structure of minimal Markov
bases and gives some characterizations. In particular, Aoki and
Takemura(\cite{aoki-takemura-2003anz}, \cite{aoki-takemura-2003metr},
\cite{aoki-takemura-2005jscs}) give the expression of the minimal Markov
bases directly (i.e., not by using algebraic algorithm) for some problems of
contingency tables. 

In applications, it is most convenient for readers to rely on algebraic
computational packages such as 4ti2 (\cite{4ti2}). 
For example, consider the problem we have seen in Section 2.1. 
For the model (\ref{eqn:Hamada-Nelder-model}), the covariate matrix
$X_0'$ is a $10\times 16$ matrix. To calculate the Markov basis for this
$X_0'$ by 4ti2, we only have to prepare a datafile
\begin{verbatim}
10 16
 1  1  1  1  1  1  1  1  1  1  1  1  1  1  1  1     
 1  1  1  1  1  1  1  1 -1 -1 -1 -1 -1 -1 -1 -1     
 1  1  1  1 -1 -1 -1 -1  1  1  1  1 -1 -1 -1 -1     
 1  1 -1 -1  1  1 -1 -1  1  1 -1 -1  1  1 -1 -1     
 1 -1  1 -1  1 -1  1 -1  1 -1  1 -1  1 -1  1 -1     
 1 -1  1 -1 -1  1 -1  1 -1  1 -1  1  1 -1  1 -1     
 1 -1 -1  1  1 -1 -1  1  1 -1 -1  1  1 -1 -1  1     
 1 -1 -1  1 -1  1  1 -1  1 -1 -1  1 -1  1  1 -1     
 1  1 -1 -1  1  1 -1 -1 -1 -1  1  1 -1 -1  1  1     
 1 -1  1 -1 -1  1 -1  1  1 -1  1 -1 -1  1 -1  1     
\end{verbatim}
and run the command {\tt markov}. Then the list of a minimal Markov
basis is instantly given as
\begin{verbatim}
35 16
-1 -1  0  0  1  1  0  0  0  0  1  1  0  0 -1 -1
-1 -1  0  1  1  1 -1  0  0  0  1  0  0  0  0 -1
-1 -1  1  0  1  1  0 -1  0  0  0  1  0  0 -1  0
....
\end{verbatim}
The above output shows that a  minimal Markov bases for this
$X'_0$ consists of $35$ moves, which corresponds to each row.

Using a minimal Markov basis, 
we perform the likelihood ratio test based on 
(\ref{eq:g2}) in Section 2.1.
After $100,000$ burn-in steps, we construct $1,000,000$ Monte Carlo
samples. In contrast to the asymptotic $p$ value $0.0040$, 
the estimated $p$ value is $0.032$, with estimated standard
deviation $0.0045$, where we use a batching method to obtain an estimate
of variance, see \cite{hastings-1970biometrika} and \cite{ripley-1987}.
Figure 1 shows a histogram of the Monte Carlo sampling generated from the
exact conditional distribution of the likelihood ratio statistic under
the null hypothesis, along with the corresponding asymptotic
distribution $\chi^2_6$.
\begin{figure*}[htbp]
\begin{center}
\includegraphics[width=11cm, height=7cm]{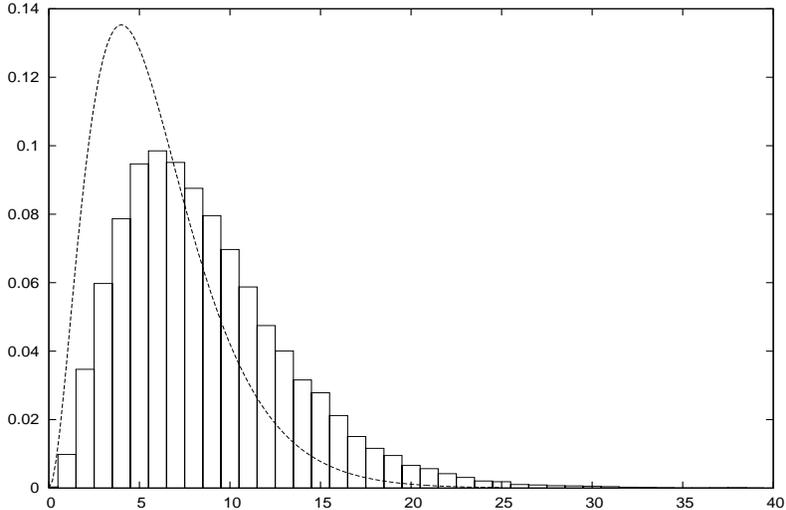}
\caption{Asymptotic and Monte Carlo estimated exact distribution}
\end{center}
\end{figure*}
Figure 1 shows that the asymptotic distribution understates the
probability that the values of the test statistic for the samples is
not less than the observed value, and overemphasizes the significance.

\subsection{Markov bases and corresponding models for $2^{p-q}$
  contingency tables}
In this section, we investigate relationships between contingency
tables and fractional factorial designs with $2^{p-q}$ runs.  As noted
in Section 1, Markov bases have been mainly investigated in the
context of contingency tables. For example, \cite{aoki-takemura-2006aism}
gives an expression of minimal Markov bases of all the hierarchical
models for $2^4$ contingency tables. This list may be sufficient in
application for the analysis of $2^4$ contingency tables, since the
hierarchical model is the most natural class of models to be
considered. We show in this section that, considering the fractional
factorial designs, we encounter some new models and Markov bases,
which do not correspond to hierarchical models of contingency tables.

\paragraph*{Fractional factorial designs with $8$ runs.}
First, we consider fractional factorial designs with $8$ runs, i.e.,
the case of $p - q = 3$. The most frequently used designs are listed in
Table 4.
\begin{table*}
\begin{center}
\caption{Eight-run $2^{p-q}$ fractional factorial designs ($p - q = 3$)}
\begin{tabular}{ccl}\hline
Number of factors $p$ & Resolution & Design Generators\\ \hline
4 & IV & $\BD = \BA\BB\BC$ \\ 
5 & III & $\BD = \BA\BB,\ \BE = \BA\BC$ \\
6 & III & $\BD = \BA\BB,\ \BE = \BA\BC,\ \BF = \BB\BC$ \\
7 & III & $\BD = \BA\BB,\ \BE = \BA\BC,\ \BF = \BB\BC, \BG =
 \BA\BB\BC$\\ \hline
\end{tabular}
\end{center}
\end{table*}
We clarify the relationships between these
designs and the models of $2^3$ contingency tables $\By = (y_{ijk}),\ 1\leq 
i,j,k\leq 2$, for the Poisson observations, and the models of $2^4$
contingency tables $\By = (y_{ijk\ell}),\ 1\leq i,j,k,\ell\leq 2$, for
the binomial observations. 

In the case of Poisson observations, we write eight observations as if
they are the frequencies of $2^3$ contingency table, i.e., 
\[
 \By = (y_{111}, y_{112}, y_{121}, y_{122}, y_{211}, y_{212}, y_{221},
 y_{222})'.
\]
In the case of $p = 5$, for example, the design and the observations are
given as follows.
\begin{center}
\begin{tabular}{ccccccc}\hline
& \multicolumn{5}{c}{Factor} & \\
Run & A & B & C & D & E & $y$ \\ \hline
1 & 1 & 1 & 1 & 1 & 1 & $y_{111}$\\
2 & 1 & 1 & 2 & 1 & 2 & $y_{112}$\\
3 & 1 & 2 & 1 & 2 & 2 & $y_{121}$\\
4 & 1 & 2 & 2 & 2 & 1 & $y_{122}$\\
5 & 2 & 1 & 1 & 2 & 1 & $y_{211}$\\
6 & 2 & 1 & 2 & 2 & 2 & $y_{212}$\\
7 & 2 & 2 & 1 & 1 & 2 & $y_{221}$\\
8 & 2 & 2 & 2 & 1 & 1 & $y_{222}$\\ \hline
\end{tabular}
\end{center}
For this type of data, we define $\nu$-dimensional parameter $\beta$ and
the covariate matrix $X$ according to an appropriate model we consider, 
as explained in Section 2. First consider the simple 
main effect model $\BA/\BB/\BC/\BD/\BE$ ($\nu = 5$). 
To test this model against various alternatives,
Markov chain Monte Carlo testing procedure needs a Markov
basis for the covariate matrix
\[
 X_0' = \left(
\begin{array}{cccccccc}
+1 & +1 & +1 & +1 & +1 & +1 & +1 & +1\\
+1 & +1 & +1 & +1 & -1 & -1 & -1 & -1\\
+1 & +1 & -1 & -1 & +1 & +1 & -1 & -1\\
+1 & -1 & +1 & -1 & +1 & -1 & +1 & -1\\
+1 & +1 & -1 & -1 & -1 & -1 & +1 & +1\\
+1 & -1 & +1 & -1 & -1 & +1 & -1 & +1
\end{array}
\right).
\]
Note that each row of the $X_0'\By$ corresponds to the sufficient
statistic under the null model $\BA/\BB/\BC/\BD/\BE$. In this case, the
sufficient statistic is given as
\begin{equation}
\begin{array}{l}
 y_{\cdot\cdot\cdot},\ y_{1\cdot\cdot},\ y_{2\cdot\cdot},\ 
  y_{\cdot 1\cdot},\ y_{\cdot 2\cdot},\  y_{\cdot\cdot 1},\
  y_{\cdot\cdot 2},\\
  y_{11\cdot}+y_{22\cdot},\ y_{12\cdot}+y_{21\cdot},\ 
y_{1\cdot 1}+y_{2\cdot 2},\ y_{1\cdot 2}+y_{2\cdot 1},
\end{array}
\label{eqn:marginals-for-ABCDE}
\end{equation}
where we use the notations such as
\[
 y_{\cdot\cdot\cdot} = \sum_{i=1}^2\sum_{j=1}^2\sum_{k=1}^2y_{ijk},\ 
 y_{i\cdot\cdot} = \sum_{j=1}^2\sum_{k=1}^2y_{ijk},\ 
 y_{ij\cdot} = \sum_{k=1}^2y_{ijk}.
\]
Here we see that the sufficient statistic
(\ref{eqn:marginals-for-ABCDE}) is 
nothing but a well-known sufficient statistic for the {\it conditional
independence model} $\BA\BB/\BA\BC$, given as
\begin{equation}
\{y_{ij\cdot}\},\ \{y_{i\cdot k}\},\ i,j,k = 1,2. 
\label{eqn:marginals-for-AB-AC}
\end{equation}
The one-to-one relation between (\ref{eqn:marginals-for-ABCDE}) and 
(\ref{eqn:marginals-for-AB-AC}) is easily shown as
\begin{equation}
\begin{array}{l}
 y_{ij\cdot} = \displaystyle\frac{y_{i\cdot\cdot} + y_{\cdot j\cdot} -
  (y_{ij^*\cdot} +
 y_{i^*j\cdot})}{2},\ 
 y_{i\cdot k} = \frac{y_{i\cdot\cdot} + y_{\cdot\cdot k} - (y_{i\cdot k^*} +
 y_{i^*\cdot k})}{2},
\end{array}
\label{eqn:correspondence}
\end{equation}
where $\{i,i^*\},\{j,j^*\}$ and $\{k,k^*\}$ are distinct indices,
respectively.
This correspondence is, of course, due to the aliasing relation 
$\BD = \BA\BB,\ \BE = \BA\BC$. 

Next consider another model. Since there 
are eight observations, we can estimate eight parameters at most (in the 
saturated model). Since the saturated model cannot be tested, 
let us 
consider the models of seven parameters, i.e., case of $\nu = 6$. 
If we restrict our attention to the hierarchical models, 
five main effects and
one of the two-factor interaction effects, $\BB\BC, \BB\BE, \BC\BD, \BD\BE$, 
can be included in the models, since the aliasing relation is given as
\begin{eqnarray*}
&&\BA = \BB\BD = \BC\BE, \ 
\BB = \BA\BD, \ 
\BC = \BA\BE, \ 
\BD = \BA\BB, \ 
\BE = \BA\BC, \\
&&\BB\BC = \BD\BE, \ 
\BB\BE = \BC\BD= \BA\BB\BC.
\end{eqnarray*}
If our null model includes $\BB\BC$ or $\BD\BE$, i.e., if our null model is written 
as $\BA/\BB\BC/\BD/\BE$ or $\BA/\BB/\BC/\BD\BE$, we add the 
column
\[
(+1\ -1\ -1\ +1\ +1\ -1\ -1\ +1)'
\]
to the covariate matrix $X_0$. In this case, the sufficient statistic
under the null model includes $y_{\cdot 11}+y_{\cdot 22}$ and $y_{\cdot
12}+y_{\cdot 21}$
in addition to (\ref{eqn:marginals-for-ABCDE}), which is nothing but 
a well-known sufficient statistic for the {\it no three-factor
interaction model}, 
$\BA\BB/\BA\BC/\BB\BC$,
\[
\{y_{ij\cdot}\},\ \{y_{i\cdot k}\},\ \{y_{\cdot jk}\},\ i,j,k = 1,2,
\]
by the similar relations to (\ref{eqn:correspondence}). 

On the other hand, if our null model includes $\BB\BE$ or $\BC\BD$, i.e., if 
our null model is written as $\BA/\BB\BE/\BC/\BD$ or $\BA/\BB/\BC\BD/\BE$, 
we have to add the column
\[
(+1\ -1\ -1\ +1\ -1\ +1\ +1\ -1)'
\]
to the covariate matrix $X_0$. In this case, the sufficient statistic
under the null model includes
$y_{111}+y_{122}+y_{212}+y_{221}$ and $y_{112}+y_{121}+y_{211}+y_{222}$ 
in addition to (\ref{eqn:marginals-for-ABCDE}).
This is one of the models which do not have corresponding models in the 
hierarchical models of three-way contingency tables. We write this new
model as $\BA\BB/\BA\BC + (\BA\BB\BC)$. The sufficient statistic for
this model is
\[
 \{y_{ij\cdot}\},\ \{y_{i\cdot k}\},\ y_{111}+y_{122}+y_{212}+
y_{221},\ y_{112}+y_{121}+y_{211}+y_{222},\ i,j,k = 1,2.
\] 

Similarly, we can specify the corresponding models of three-way
contingency tables (for the factors $\BA,\BB,\BC$) to all the possible
models for the designs of
Table 4, as if the observations are the frequencies of $2^3$ contingency
tables. The result is summarized in Table 5.

\begin{table*}
\begin{center}
\caption{Eight-run $2^{p-q}$ fractional factorial designs and the
 corresponding models of three-way contingency tables ($p - q = 3$)}
\medskip
\begin{tabular}{lll}\hline
\multicolumn{3}{l}{Design : $p = 4$, $\BD = \BA\BB\BC$} \\ 
$\nu$ & Null model & Corresponding model of $2^3$ table \\ \hline
$4$ & $\BA/\BB/\BC/\BD$ & $\BA/\BB/\BC + (\BA\BB\BC)^a$ \\ 
$5$ & $\BA\BB/\BC/\BD$ & $\BA\BB/\BC + (\BA\BB\BC)^a$ \\
$6$ & $\BA\BB/\BA\BC/\BD$ & $\BA\BB/\BA\BC + (\BA\BB\BC)^a$ \\ \hline
\multicolumn{3}{l}{Design : $p = 5$, $\BD = \BA\BB,\ \BE = \BA\BC$} \\ 
$\nu$ & Null model & Corresponding model of $2^3$ table \\ \hline
$5$ & $\BA/\BB/\BC/\BD/\BE$ & $\BA\BB/\BA\BC$ \\
$6$ & $\BA/\BB\BC/\BD/\BE$ & $\BA\BB/\BA\BC/\BB\BC$ \\ 
    & $\BA/\BB\BE/\BC/\BD$ & $\BA\BB/\BA\BC + (\BA\BB\BC)^a$ \\ \hline
\multicolumn{3}{l}{Design : $p = 6$, $\BD = \BA\BB,\ \BE = \BA\BC,\ \BF
 = \BB\BC$} \\ 
$\nu$ & Null model & Corresponding model of $2^3$ table \\ \hline
$6$ & $\BA/\BB/\BC/\BD/\BE/\BF$ & $\BA\BB/\BA\BC/\BB\BC$ \\ \hline
\end{tabular}
\ \\
\medskip
${}^a$The sufficient statistic for $(\BA\BB\BC)$ is $y_{111}+y_{122}+y_{212}+
y_{221},\ y_{112}+y_{121}+y_{211}+y_{222}$.
\end{center}
\end{table*}

\bigskip

In the case of Binomial observations, there are $16$ observations.
Similarly to the Poisson case, 
we treat the observations as if they are the frequencies of $2^4$
contingency table. In the case of $p = 5$, for example, the design and the
observations are given as follows.
\begin{center}
\begin{tabular}{cccccccc}\hline
& \multicolumn{5}{c}{Factor} & & \\
Run & A & B & C & D & E & \multicolumn{2}{c}{$y$} \\ \hline
1 & 1 & 1 & 1 & 1 & 1 & $y_{1111}$ & $y_{1112}$\\
2 & 1 & 1 & 2 & 1 & 2 & $y_{1121}$ & $y_{1122}$\\
3 & 1 & 2 & 1 & 2 & 2 & $y_{1211}$ & $y_{1212}$\\
4 & 1 & 2 & 2 & 2 & 1 & $y_{1221}$ & $y_{1222}$\\
5 & 2 & 1 & 1 & 2 & 1 & $y_{2111}$ & $y_{2112}$\\
6 & 2 & 1 & 2 & 2 & 2 & $y_{2121}$ & $y_{2122}$\\
7 & 2 & 2 & 1 & 1 & 2 & $y_{2211}$ & $y_{2212}$\\
8 & 2 & 2 & 2 & 1 & 1 & $y_{2221}$ & $y_{2222}$\\ \hline
\end{tabular}
\end{center}
For this type of data, we also specify parameter
$\beta$ and the covariate matrix according to the appropriate models, 
by replacing $X$ by $\widetilde{X}$ of (\ref{eqn:lawlence-lifting}).
Note that the elements of $\By$ is ordered as
\[
\By = (y_{1111},y_{1121},\ldots,y_{2211},y_{2221},y_{1112},y_{1122},
\ldots,y_{2212},y_{2222})' .
\]
Accordingly, correspondences to the models of $2^4$
contingency tables are easily obtained and the result is given in Table 6. 

Table 6 is automatically converted from Table 5 as follows. 
By definition, $\BD$ is added to all the generating sets.
Note also that the sufficient statistic for each model 
includes $\{y_{ijk\cdot}\}, 1\leq i,j,k \leq 2$, by definition, which
yields Table 6.
Therefore the models which do not include all of $\BA\BB, \BA\BC$ and
$\BB\BC$ do not correspond to hierarchical models. 

\begin{table*}
\begin{center}
\caption{Eight-run $2^{p-q}$ fractional factorial designs and the
 corresponding models of three-way contingency tables ($p - q = 3$)}
\medskip
\begin{tabular}{lll}\hline
\multicolumn{3}{l}{Design : $p = 4$, $\BD = \BA\BB\BC$} \\ 
$\nu$ & Null model & Corresponding model of $2^4$ table \\ \hline
$4$ & $\BA/\BB/\BC/\BD$ & $\BA\BD/\BB\BD/\BC\BD + (\BA\BB\BC)^a +
 (\BA\BB\BC\BD)^b$ \\ 
$5$ & $\BA\BB/\BC/\BD$ & $\BA\BB\BD/\BC\BD + (\BA\BB\BC)^a +
 (\BA\BB\BC\BD)^b$ \\
$6$ & $\BA\BB/\BA\BC/\BD$ & $\BA\BB\BD/\BA\BC\BD + (\BA\BB\BC)^a +
 (\BA\BB\BC\BD)^b$ \\ \hline
\multicolumn{3}{l}{Design : $p = 5$, $\BD = \BA\BB,\ \BE = \BA\BC$} \\ 
$\nu$ & Null model & Corresponding model of $2^4$ table \\ \hline
$5$ & $\BA/\BB/\BC/\BD/\BE$ & $\BA\BB\BD/\BA\BC\BD + (\BA\BB\BC)^a$  \\
$6$ & $\BA/\BB\BC/\BD/\BE$ & $\BA\BB\BD/\BA\BC\BD/\BB\BC\BD/\BA\BB\BC$ \\ 
    & $\BA/\BB\BE/\BC/\BD$ & $\BA\BB\BD/\BA\BC\BD + (\BA\BB\BC)^a +
 (\BA\BB\BC\BD)^b$ \\ \hline
\multicolumn{3}{l}{Design : $p = 6$, $\BD = \BA\BB,\ \BE = \BA\BC,\ \BF
 = \BB\BC$} \\ 
$\nu$ & Null model & Corresponding model of $2^4$ table \\ \hline
$6$ & $\BA/\BB/\BC/\BD/\BE/\BF$ &
 $\BA\BB\BD/\BA\BC\BD/\BB\BC\BD/\BA\BB\BC$ \\ \hline
\end{tabular}
\end{center}
${}^a$The sufficient statistic for $(\BA\BB\BC)$ is $\{y_{ijk\cdot}\},\
 i,j,k = 1,2.$\\
${}^b$The sufficient statistic for $(\BA\BB\BC\BD)$ is
 $y_{111\ell}+y_{122\ell}+y_{212\ell}+ y_{221\ell},\\
y_{112\ell}+y_{121\ell}+y_{211\ell}+y_{222\ell},\ \ell = 1,2$.
\end{table*}

\paragraph*{Fractional factorial designs with $16$ runs.}
Next we consider fractional factorial designs with $16$ runs, i.e., the
case of $p - q = 4$. Table 7 is a list of sixteen-run $2^{p-q}$
fractional factorial designs ($p - q = 4, p \leq 10$) from Section 4
of \cite{wu-hamada-2000}. 
\begin{table*}
\begin{center}
\caption{Sixteen-run $2^{p-q}$ fractional factorial designs ($p - q = 4$)}
\medskip
\begin{tabular}{ccl}\hline
Number of factors $p$ & Resolution & Design Generators\\ \hline
5 & V   & $\BE = \BA\BB\BC\BD$ \\
6 & IV  & $\BE = \BA\BB\BC, \BF = \BA\BB\BD$\\
7 & IV  & $\BE = \BA\BB\BC, \BF = \BA\BB\BD, \BG = \BA\BC\BD$\\
8 & IV  & $\BE = \BA\BB\BC, \BF = \BA\BB\BD, \BG = \BA\BC\BD$\\
  &     & $\BH = \BB\BC\BD$\\
9 & III & $\BE = \BA\BB\BC, \BF = \BA\BB\BD, \BG = \BA\BC\BD$\\
  &     & $\BH = \BB\BC\BD, \BJ = \BA\BB\BC\BD$\\
10 & III & $\BE = \BA\BB\BC, \BF = \BA\BB\BD, \BG = \BA\BC\BD$\\
   &     & $\BH = \BB\BC\BD, \BJ = \BA\BB\BC\BD, \BK = \BC\BD$\\ \hline
\end{tabular}
\end{center}
\end{table*}
By the similar considerations to the $8$ run cases,
we can seek the corresponding models of 
$2^4$ contingency tables for
the Poisson observations, and models of $2^5$ contingency tables for
the Binomial observations. Since the modeling for Binomial observations
can be easily obtained from the Poisson case as we have seen, we only
consider the Poisson case here. 

Since at most sixteen parameters are estimable for the sixteen-run
designs, we can consider various models of main effects and interaction
effects. For example, the saturated model of the $p
= 5$ design, $\BE = \BA\BB\BC\BD$, can include all the main and two-factor
interactions,
\[
 \BA\BB/\BA\BC/\BA\BD/\BA\BE/\BB\BC/\BB\BD/\BB\BE/\BC\BD/\BC\BE/\BD\BE.
\]
Note that for the models of $p = 5,6,7,8$ in Table 7, each main effect and 
two-factor interaction is estimable. (On the other hand, for the resolution III 
models of $p = 9,10$, some of the two-factor interactions are 
not estimable.) Among the models which include all the main effects and some 
of the two-factor interaction effects, some models have the corresponding hierarchical model
in the $2^4$ contingency tables if we write the sixteen observations as 
$\By = \{y_{ijk\ell}\}, i,j,k,\ell = 1,2$. For example, for the $p = 6$ design 
of $\BE = \BA\BB\BC, \BF = \BA\BB\BD$, the model of $6$ main effects and $5$ 
two-factor interaction effects, 
\[
\BA\BB/\BA\BC/\BA\BD/\BB\BC/\BB\BD/\BE/\BF,
\]
has a corresponding model of $\BA\BB\BC/\BA\BB\BD$
in the $2^4$ contingency tables. By the aliasing relation
\begin{eqnarray*}
&&\BA\BB = \BC\BE = \BD\BF, \ 
\BA\BC = \BB\BE, \ 
\BA\BD = \BB\BF, \ 
\BA\BE = \BB\BC, \\
&&\BA\BF = \BB\BD, \ 
\BC\BD = \BE\BF, \
\BC\BF = \BD\BE = \BA\BB\BC\BD,
\end{eqnarray*}
it is seen that there are $3\cdot 2\cdot 2\cdot 2\cdot 2 = 48$ distinct models such
as
\[\begin{array}{l}
\BA\BB/\BA\BC/\BA\BD/\BA\BE/\BA\BF/\BE/\BF,\\
\BA\BB/\BA\BC/\BA\BD/\BA\BE/\BB\BD/\BE/\BF,\\
\BA\BB/\BA\BC/\BA\BD/\BB\BC/\BA\BF/\BE/\BF,\\
\BA\BB/\BA\BC/\BA\BD/\BB\BC/\BB\BD/\BE/\BF,\\
\multicolumn{1}{c}{\vdots}\\
\BD\BF/\BB\BE/\BB\BF/\BB\BC/\BA\BF/\BE/\BF,\\
\BD\BF/\BB\BE/\BB\BF/\BB\BC/\BB\BD/\BE/\BF,
\end{array}
\]
which correspond to the model of $\BA\BB\BC/\BA\BB\BD$
in the $2^4$ contingency tables. 
By the similar considerations, we can specify all the models for the
designs of Table 7
which correspond some hierarchical models in the $2^4$ contingency
tables. The result is shown in Table 8.
\begin{table*}
\begin{center}
\caption{Sixteen-run $2^{p-q}$ fractional factorial designs and the
 corresponding hierarchical models of $2^4$ contingency tables ($p - q = 4$)}
\medskip
\begin{tabular}{llll}\hline
\multicolumn{4}{l}{Design : $p = 6$, $\BE = \BA\BB\BC, \BF = \BA\BB\BD$} \\ 
$\nu$ & Representative & Num. of the     & Corresponding hierarchical\\ 
      & null model     & null models     & model of $2^4$ table \\ \hline
$11$ & $\BA\BB/\BA\BC/\BA\BD/\BB\BC/\BB\BD/\BE/\BF$ & $48$ &
 $\BA\BB\BC/\BA\BB\BD$\\ 
$12$ & $\BA\BB/\BA\BC/\BA\BD/\BB\BC/\BB\BD/\BC\BD/\BE/\BF$ & $96$ & 
$\BA\BB\BC/\BA\BB\BD/\BC\BD$\\ \hline
\multicolumn{4}{l}{Design : $p = 7$, $\BE = \BA\BB\BC, \BF = \BA\BB\BD, \BG = \BA\BC\BD$} \\ 
$\nu$ & Representative & Num. of the     & Corresponding hierarchical\\ 
      & null model     & null models     & model of $2^4$ table \\ \hline
$11$ & $\BA\BB/\BA\BC/\BA\BD/\BB\BC/\BB\BD/\BC\BD/\BE/\BF/\BG$ & $3^6 = 729$ & 
$\BA\BB\BC/\BA\BB\BD/\BA\BC\BD$\\ \hline
\multicolumn{4}{l}{Design : $p = 8$, $\BE = \BA\BB\BC, \BF = \BA\BB\BD,
 \BG = \BA\BC\BD, 
\BH = \BB\BC\BD$} \\ 
$\nu$ & Representative & Num. of the     & Corresponding hierarchical\\ 
      & null model     & null models     & model of $2^4$ table \\ \hline
$11$ & $\BA\BB/\BA\BC/\BA\BD/\BB\BC/\BB\BD/\BC\BD/\BE/\BF/\BG$ & $4^6 =
 4096$ & 
$\BA\BB\BC/\BA\BB\BD/\BA\BC\BD/\BB\BC\BD$\\ \hline
\end{tabular}
\end{center}
\end{table*}

One of the merits to specify corresponding hierarchical models of
contingency tables is a possibility to make use of general already
known results on Markov bases of contingency tables. 
For example, \cite{dobra-2003bernoulli} shows that a Markov basis can be
constructed by 
the primitive moves, i.e., degree 2 moves, for the decomposable
graphical models in the contingency tables.
In our designed experiments, therefore, Markov basis for the models
which correspond to 
decomposable graphical models of contingency tables can be
constructed only by the 
primitive moves. 
Among the results of Table 5,6 and 8, there are two models which
corresponds to decomposable graphical models in the contingency tables.
We can confirm that minimal Markov bases for these models consist
of primitive moves as follows.  We use the following  notational
convention for a move $\Bz$.  Consider $2^3$ case
$\Bz=(z_{ijk})$.   If $z_{i_1 j_1 k_1}=z_{i_2 j_2 k_2}=+1$, 
$z_{i_3 j_3 k_3}=z_{i_4 j_4 k_4}=+1$, and other elements are zeros,
then we denote $\Bz$ as
\[
(i_1 j_1 k_1)(i_2 j_2 k_2) - (i_3 j_3 k_3)(i_4 j_4 k_4).
\]
Similar notation is used for $2^4$ case.

\begin{itemize}
\item $2^{5-2}$ fractional factorial design of $\BD = \BA\BB, \BE =
      \BA\BC$:\\
      The main effects model $\BA/\BB/\BC/\BD/\BE$ corresponds to the
      decomposable graphical model $\BA\BB/\BA\BC$ of the $2^3$
      contingency tables. 
      This is a conditional independence model 
      between $\BB$ and $\BC$ given $\BA$ and a minimal Markov basis is
      constructed by primitive moves as
\[
(111)(122)-(112)(121),\ (211)(222)-(212)(221).
\]
\item $2^{6-2}$ fractional factorial design of $\BE = \BA\BB\BC, \BF =
      \BA\BB\BD$:\\
      The model $\BA\BB/\BA\BC/\BA\BD/\BB\BC/\BB\BD/\BE/\BF$ corresponds
      to the decomposable graphical model $\BA\BB\BC/\BA\BB\BD$ of the 
      $2^4$ contingency tables.
      This is a conditional independence model 
      between $\BC$ and $\BD$ given $\{\BA,\BB\}$ and a minimal Markov
      basis is again 
      constructed by primitive moves as
\[\begin{array}{c}
(1111)(1122)-(1112)(1121),\ (1211)(1222)-(1212)(1221),\\
(2111)(2122)-(2112)(2121),\ (2211)(2222)-(2212)(2221).
\end{array}
\]

\end{itemize}

For the other designs of Table 7 ($p = 5,9,10$), all the models include
the sufficient statistic
\[\begin{array}{c}
y_{1111}+y_{1122}+y_{1212}+y_{1221}+y_{2112}+y_{2121}+y_{2211}+y_{2222},\\
y_{1112}+y_{1121}+y_{1211}+y_{1222}+y_{2111}+y_{2122}+y_{2212}+y_{2221},
\end{array}
\]
and therefore have no corresponding hierarchical models in the $2^4$
contingency tables.
For example, the sufficient statistic of the the main effect models for 
$2^{5-1}$ design of $\BE = \BA\BB\BC\BD$ is
\[\begin{array}{l}
\{y_{i\cdot\cdot\cdot}\},\ \{y_{\cdot j\cdot\cdot}\},\ \{y_{\cdot\cdot
 k\cdot}\},\ 
\{y_{\cdot\cdot\cdot \ell}\},\ i,j,k,\ell = 1,2,\\ 
y_{1111}+y_{1122}+y_{1212}+y_{1221}+y_{2112}+y_{2121}+y_{2211}+y_{2222},\\
y_{1112}+y_{1121}+y_{1211}+y_{1222}+y_{2111}+y_{2122}+y_{2212}+y_{2221},
\end{array}
\]
and the sufficient statistic of the main effect models for 
$2^{10-1}$ design of 
\[
\BE = \BA\BB\BC,\ \BF = \BA\BB\BD,\ \BG = \BA\BC\BD,\ \BH = \BB\BC\BD,\ 
\BJ = \BA\BB\BC\BD,\ \BK = \BC\BD
\]
is
\[\begin{array}{l}
\{y_{ijk\cdot}\},\ \{y_{ij\cdot\ell}\},\ \{y_{i\cdot k\ell}\},\ 
\{y_{\cdot jk\ell}\},\ i,j,k,\ell = 1,2,\\ 
y_{1111}+y_{1122}+y_{1212}+y_{1221}+y_{2112}+y_{2121}+y_{2211}+y_{2222},
\\
y_{1112}+y_{1121}+y_{1211}+y_{1222}+y_{2111}+y_{2122}+y_{2212}+y_{2221}
.
\end{array}
\]

\section{Discussion}
In this paper, we consider Markov chain Monte Carlo tests for the
factor effects in the designed experiments. As is noted in Section 1,
Markov chain Monte Carlo procedure is a valuable tool when the
adequacy of traditional large-sample tests is doubtful and 
the enumeration of the conditional sample space is infeasible.  
Since a closed form expression of the null
distribution for the conditional tests considered in Section 2 is not
available in general, Markov chain Monte Carlo procedure is valuable
in the settings of this paper.  Computational experience given in
Section 3.1 shows efficacy of our method.

To perform Markov chain Monte Carlo tests, it is often problematic to
calculate a Markov basis. Current algorithms may take a very long time
to compute Markov basis for 64 run or larger experiments.  For the
designs of $16$ or $32$ runs (for the Poisson models), a software such
as 4ti2 works quite well and very practical.  Nevertheless, the
arguments and theoretical considerations in Section 3.2 seem
important.  One of the merits to specify the corresponding models of
$2^p$ contingency tables is a possibility to make use of general
results for the Markov bases of contingency tables as shown in Section
3.2.

It is also important to consider more complicated designs and give
appropriate Markov bases for them, such as designs with three levels
or balanced incomplete block designs.

The designed experiment is one of the areas in statistics where 
the applications of the theory of Gr\"obner basis are first considered. 
See the works \cite{pistone-wynn-1996biometrika}, 
\cite{robbiano-rogantin-1998} and \cite{pistone-riccomagno-wynn-2000}. 
In these works, the design is represented as the variety
for the set of polynomial equations. On the other hand, this
manuscript gives another application of Gr\"obner basis theory
to the designed experiments
by considering Markov chain Monte Carlo tests for a discrete response variable.

\bibliographystyle{plain}
\bibliography{experiment}

\end{document}